# Localized Curvature Domination and Rigidity of Harmonic Maps

Sergey E. Stepanov

Department of Mathematics, Finance University
125468 Moscow, Leningradsky Prospect, 49-55, Russian Federation
e-mail: *s.e.stepanov@mail.ru*

**Abstract.** We establish a localized Bochner-type rigidity theorem for harmonic maps between Riemannian manifolds. Let $f: (M, g) \to (\bar{M}, \bar{g})$ be a harmonic map from a compact manifold. Instead of assuming a global nonpositivity condition on the sectional curvature of the target, we impose a curvature bound localized to the image $f(M)$, expressed via the maximal sectional curvature encountered along the image.
We prove that if the minimal Ricci curvature of the domain dominates this image-dependent curvature bound in a sharp quantitative pinching inequality involving the maximal energy density of $f$, then the map is constant. At the critical threshold, we obtain a homothetic classification: the differential is parallel and the image is totally geodesic.
The result replaces global curvature sign assumptions with an image-dependent curvature domination principle and yields a localized analogue of Yano–Ishihara-type rigidity.

## 1. Introduction

Harmonic maps constitute one of the most prolific fields of geometric analysis, situated at the intersection of differential geometry, the theory of partial differential equations, and theoretical physics (see, for example, [1]; [2] and [3]).

Harmonic maps between Riemannian manifolds $f: (M, g) \to (\bar{M}, \bar{g})$ arise as critical points of the energy functional (see [4])

$$E(f) = \frac{1}{2} \int_M |df|^2 \, d\mathrm{vol}_g,$$

and their rigidity properties are governed by a subtle interaction between the curvature of the domain, the curvature of the target, and global geometric constraints. Since the seminal work of Eells–Sampson (see [4]), curvature sign conditions have played a decisive role in rigidity theory. In the compact setting, Eells–Sampson showed that if the target manifold has nonpositive sectional curvature, every homotopy class contains a harmonic representative obtained via the heat flow, and strict negativity of sectional curvature enforces strong rigidity within the homotopy class. In a complementary direction, Yau and Schoen established a Liouville-type rigidity theorem (see [5]): if the domain has nonnegative Ricci curvature and the compact target has nonpositive sectional curvature, then every harmonic map must be constant. In both cases the non-positivity of sectional curvature of the target is the fundamental geometric assumption driving the vanishing phenomenon.

The analytic core of these results is the Bochner–Weitzenböck identity (see [4])

$$\Delta e(f) = |\nabla df|^2 + \langle \text{Ric}_M(df), df \rangle - \langle \bar{R}(df, df)df, df \rangle,$$

where $e(f) = \frac{1}{2}|df|^2$. When curvature sign conditions render the right-hand side nonnegative, the energy density becomes subharmonic, and compactness or growth arguments force it to be constant, leading to rigidity. Classical theory therefore relies on global qualitative curvature assumptions, typically requiring $\overline{\text{Sec}} \leq 0$ on the entire target manifold (see [4 - 6]).

The purpose of this paper is to refine this rigidity mechanism by localizing curvature control to the image of the map. If $M$ is compact, then $f(M) \subset \bar{M}$ is compact, and one may define the maximal sectional curvature encountered along the image,

$$\overline{\text{Sec}}_{\text{nax}}(f(M)) = \max_{\substack{q \in f(M) \\ \sigma \subset T_q\bar{M}}} \overline{\text{Sec}}(\sigma).$$

Instead of imposing a global sign condition on $\overline{\text{Sec}}$, we assume a quantitative curvature domination inequalities of the form

$$\text{Ric}_{\min}(M) > \frac{n-1}{n} \overline{\text{Sec}}_{\max}(f(M)) \cdot \sup_M |df|^2,$$

where $\text{Ric}_{\min}(M)$ denotes the minimal Ricci curvature of the domain and $\overline{\text{Sec}}(\sigma) \geq 0$ for an arbitrary $\sigma \subset T_q\bar{M}$ at each point $p \in f(M)$. Under this assumption we prove that a harmonic map must be constant. Thus rigidity is governed not by the global geometry of the target, but by the curvature effectively experienced along the image together with a quantitative balance between intrinsic and ambient curvature contributions.

This localization principle becomes significant when the target manifold does not satisfy a global curvature sign condition. For instance, if $\bar{M}$ contains regions of positive sectional curvature but the image of the harmonic map lies entirely in a region where sectional curvature is bounded above by a constant $\kappa$, then rigidity follows provided the Ricci curvature of the domain dominates $\kappa$ in the quantitative sense above. The geometry of regions disjoint from the image plays no role in the argument. In this way, the present result extends the scope of classical Bochner-type rigidity beyond the traditional framework of globally nonpositive curvature.

A further structural feature of the theory is the analysis of the threshold case. When the curvature domination inequality becomes an equality, the Bochner identity forces the differential $df$ to be parallel, and the algebraic equalities underlying the curvature estimate imply coincidence of the singular values of $df$. Consequently, the map is either constant or homothetic, and its image is totally geodesic. This phenomenon parallels the classification result of Yano and Ishihara, (see [6]) who showed that under a global curvature pinching condition involving a uniform upper bound on sectional curvature of the target, a harmonic map from a compact manifold with positive Ricci curvature must be either constant or homothetic. The present framework refines that perspective by replacing the global curvature bound with an image-dependent one and by expressing the pinching condition in terms of the maximal energy density of the map.

Taken together, these results establish a localized curvature domination principle: quantitative control of ambient sectional curvature along the image, dominated by intrinsic Ricci curvature of the domain, yields either collapse to a constant map in the strict regime or homothetic rigidity at the critical threshold. In this sense, the classical Bochner method is preserved in analytic structure while its geometric hypotheses are significantly relaxed and localized.

In the present paper we continue the study which we began in [7] where we proved some vanishing theorems for harmonic mappings into non-negatively curved manifolds and studied their applications.

## 2. Main theorem

Let $(M, g)$ be a compact Riemannian manifold and $(\bar{M}, \bar{g})$ an arbitrary Riemannian manifold, and let

$$f: (M, g) \to (\bar{M}, \bar{g})$$

be a differentiable map. We begin by recalling a standard compactness consequence for curvature on $M$. Since the Ricci curvature tensor $\text{Ric}_M$ depends smoothly on the metric, its eigenvalues vary continuously on $M$. Because $M$ is compact, the Ricci curvature admits a global lower bound; that is, there exists a constant

$$\text{Ric}_{\min}(M) = \min_{p \in M} \min_{|X|=1} (\text{Ric}(X, X))$$

This minimum is attained at some point of $M$.

We now turn to the geometry of the target along the image of $f$. Since every differentiable map between smooth manifolds is continuous, and continuous maps send compact sets to compact sets, the image $f(M) \subset \bar{M}$ is a compact subset of $\bar{M}$. No additional assumptions on $\bar{M}$ are required for this conclusion.

To discuss sectional curvature along the image, consider the Grassmann bundle of two-dimensional planes over $\bar{M}$, $\text{Gr}_2(T\bar{M}) \to \bar{M}$. The sectional curvature function

$$\overline{\text{Sec}} : \text{Gr}_2(T\bar{M}) \to \mathbb{R}$$

is continuous, as it depends smoothly on the metric tensor and the curvature tensor of $\bar{M}$. Restricting this bundle to the compact set $f(M)$ yields $\text{Gr}_2(T\bar{M})|_{f(M)}$, which is compact because each fibre $\text{Gr}_2(T_p\bar{M})$ is compact and the base $f(M)$ is compact. Consequently, by the extreme value theorem, the sectional curvature function attains its minimum on this restricted bundle. Hence the quantity

$$\overline{\text{Sec}}_{\max}(f(M)) := \max_{\substack{p \in f(M) \\ \sigma \subset T_p \bar{M}}} \overline{\text{Sec}}(\sigma)$$

is well defined and achieved for some two-plane $\sigma \subset T_p \bar{M}$ with $p \in f(M)$.

It is important to stress that the minimization is taken over all two-dimensional planes in $T_p \bar{M}$ at points of the image, not merely over planes tangent to $f(M)$. The latter would describe curvature properties intrinsic to a submanifold, whereas the present quantity reflects the ambient geometry of $\bar{M}$ along the image of the map.

In summary, compactness of $M$ guarantees the existence of a global lower bound for its Ricci curvature and ensures that the image $f(M)$ is compact in $\bar{M}$. This compactness, in turn, implies the existence of a global maximum of the sectional curvature of $\bar{M}$ over all two-planes based at points of $f(M)$. Then the following theorem is true.

**Main Theorem** (*Rigidity under Curvature Pinching*)

*Let $f: (M, g) \to (\bar{M}, \bar{g})$ be a smooth harmonic map from an $n$-dimension connected, compact, boundaryless Riemannian manifold $(M, g)$ into a connected Riemannian manifold $(\bar{M}, \bar{g})$. Denote by*

$$\text{Ric}_{\min}(M) = \min_{p \in M} \min_{|X|=1} (\text{Ric}(X, X))$$

*the global minimum of the Ricci curvature of $(M, g)$, and define by*

$$\overline{\text{Sec}}_{\max}(f(M)) = \max_{\substack{p \in f(M) \\ \sigma \subset T_p \bar{M}}} \overline{\text{Sec}}(\sigma),$$

*the global maximum of the sectional curvature of $(\bar{M}, \bar{g})$, which is well defined by compactness of $f(M)$. Let $e(f) = \frac{1}{2} |df|^2$ denote the energy density of $f$, and set*

$$e_{\max} = \sup_M e(f).$$

*Assume that $\overline{\text{Sec}}(\sigma) \geq 0$ for an arbitrary two-plane $\sigma \subset T_p \bar{M}$ at each $p \in f(M)$ and*

$$\text{Ric}_{\min}(M) > \frac{n-1}{n} e_{\max} \cdot \overline{\text{Sec}}_{\max}(f(M)).$$

*Then $f$ is constant.*

## 3. Assessment of the Main Theorem

The theorem establishes a Bochner-type rigidity result for harmonic maps under a quantitative curvature pinching condition. Its distinguishing feature is the replacement of the classical global sign assumption on the sectional curvature of the target manifold by a localized upper bound taken only along the image $f(M)$. In particular, the curvature constant involved in the estimate is not a global invariant of the target manifold but rather the maximal sectional curvature evaluated over two-planes based at points of the image. This represents a conceptual refinement of standard rigidity frameworks.

In contrast to the classical results of Eells–Sampson (see [4]) and Yau- Schoen (see [5]), where nonpositive sectional curvature of the entire target manifold plays a decisive structural role, the present theorem permits the target to contain regions of positive sectional curvature. The rigidity conclusion is obtained through a quantitative balance inequality between the minimal Ricci curvature of the domain and the maximal sectional curvature of the target restricted to the image. The mechanism is thus intrinsically competitive: the intrinsic curvature of the domain must dominate the ambient curvature experienced by the map along its image.

Structurally, the result belongs to the class of *Bochner-type vanishing theorems*. The compactness of the domain guarantees both the attainment of the minimal Ricci curvature and the compactness of the image, which in turn ensures that the maximal sectional curvature along $f(M)$ is well defined. The key inequality

$$\text{Ric}_{\min}(M) > \frac{n-1}{n} \overline{\text{Sec}}_{\max}(f(M)) \cdot \sup_{M} \mid df \mid^2$$

expresses a curvature domination principle: positive Ricci curvature of the domain must control the potential curvature amplification effect arising from the target geometry, weighted by the maximal size of the differential of the map.

From a mathematical standpoint, the strength of the theorem lies in two principal aspects:

1. It does not impose global curvature restrictions on the target manifold.

2. It replaces qualitative sign conditions with a quantitative rigidity threshold derived from a sharp algebraic estimate in the curvature term of the Bochner formula.

At the same time, the hypothesis remains restrictive. The presence of the factor $\sup_M |df|^2$ (or equivalently $\sup_M e(f)$) makes the rigidity condition dependent on the analytic size of the map. In the particular case where the sectional curvature of the target is nonpositive along the image, the inequality reduces to a classical vanishing condition under positive Ricci curvature of the domain. Thus, the genuinely new content of the theorem emerges in situations where the target manifold admits positive sectional curvature, but its effect is quantitatively dominated by the curvature of the domain.

The coefficient $\frac{n-1}{n}$ is not incidental; it arises from the optimal algebraic inequality relating $\sum_{i<j} \lambda_i \lambda_j$ to $(\sum_i \lambda_i)^2$ in the diagonalization of the pullback metric. If this constant is indeed sharp, the theorem possesses structural optimality within the Bochner framework.

In summary, the theorem provides a geometrically coherent and technically precise rigidity statement. Its novelty lies in the localization of sectional curvature control to the image of the map and in the formulation of a quantitative curvature balance principle. The broader impact of the result will depend on the sharpness of the pinching constant and on the availability of geometrically meaningful examples situated near the critical threshold.

### 3. Proof of the Main Theorem

1. *Notation and curvature bounds*

Let $f: (M, g) \to (\bar{M}, \bar{g})$ be a smooth *harmonic* map, where $(M, g)$ is connected, compact, and without boundary, and $(\bar{M}, \bar{g})$ is connected. Define as above:

$$\text{Ric}_{\min} := \min_{p \in M} \min_{\substack{X \in T_p M \\ |X|=1}} \text{Ric}_M(X, X), \qquad \overline{\text{Sec}}_{\max} := \max_{\substack{q \in f(M) \\ \sigma \subset T_q \bar{M}}} \overline{\text{Sec}}(\sigma).$$

Since $M$ is compact, $\text{Ric}_{\inf}$ is well-defined and attained. Since $f(M)$ is compact, $\overline{\text{Sec}}_{\max}$ is well-defined and attained.

We use the convention

$$e(f):=\frac{1}{2}|df|^2, \quad e_{\max}:=\sup_M e(f), \quad \text{so that} \quad |df|^2 \leq 2e_{\max} \text{ on } M.$$

## 2. Bochner identity

For a harmonic map $f$, the Bochner–Weitzenböck formula reads (see [4, p 123] and [6, p. 506])

$$\frac{1}{2}\Delta|df|^2 = \|\nabla df\|^2 + Q(f),$$

where

$$Q(f) = \sum_{i=1}^n \text{Ric}_M(e_i, e_i)\,|df(e_i)|^2 - \sum_{i,j=1}^n \langle \bar{R}(df(e_i), df(e_j))\,df(e_j), df(e_i)\rangle$$

for any local $g$-orthonormal frame $\{e_1, \ldots, e_n\}$.

Integrating over compact $M$ (no boundary) yields (see [4]; [5] and [6])

$$\int_M (\|\nabla df\|^2 + Q(f))\,dv_g = 0. \tag{1}$$

## 3. Diagonalization of the pullback metric.

Fix a point $x \in M$. Choose a local $g$-orthonormal frame $\{e_1, \ldots, e_n\}$ around $x$ such that the pullback metric

$$g^* := f^*\bar{g}$$

is diagonal at $x$, i.e. $\bar{g}(df(e_i), df(e_j)) = \lambda_i\,\delta_{ij}$, $\lambda_i \geq 0$. Set $S = S(x)$ for

$$S := \sum_{i=1}^n \lambda_i = \text{tr}_g(g^*) = |df|^2.$$

## 4. Lower bound for the Ricci term.

By definition of $\text{Ric}_{\inf}$, we have (see [6])

$$\sum_{i=1}^n \text{Ric}(e_i, e_i)\,|df(e_i)|^2 = \sum_{i=1}^n \lambda_i\,\text{Ric}(e_i, e_i) \geq \text{Ric}_{\min} \sum_{i=1}^n \lambda_i = \text{Ric}_{\min} \cdot S.$$

## 5. Upper bound for the curvature term of the target

Let $u_i := df(e_i)$. Since $\bar{g}(u_i, u_j) = 0$ for $i \neq j$ and $|u_i|^2 = \lambda_i$, we have for $i \neq j$

$$\langle \bar{R}(u_i, u_j) u_j, u_i \rangle = \overline{\text{Sec}}(u_i, u_j) \ |u_i|^2 |u_j|^2 = \overline{\text{Sec}}(u_i, u_j) \lambda_i \lambda_j.$$

Moreover, the term with $i = j$ vanishes. Hence (see also [6])

$$\sum_{i,j=1}^{n} \langle \bar{R}(u_i, u_j) u_j, u_i \rangle = 2 \sum_{1 \leq i < j \leq n} \sec(u_i, u_j) \lambda_i \lambda_j \leq 2 \overline{\text{Sec}}_{\max} \sum_{i<j} \lambda_i \lambda_j.$$

Using the identity

$$\sum_{i<j} \lambda_i \lambda_j = \frac{S^2 - \sum_{i=1}^{n} \lambda_i^2}{2} \leq \frac{S^2 - \frac{1}{n}S^2}{2} = \frac{n-1}{2n} S^2,$$

where the inequality follows from Cauchy–Schwarz $\sum \lambda_i^2 \geq \frac{1}{n}(\sum \lambda_i)^2$, we obtain

$$\sum_{i,j=1}^{n} \langle \bar{R}(u_i, u_j) u_j, u_i \rangle \leq 2 \overline{\text{Sec}}_{\max} \cdot \frac{n-1}{2n} S^2 = \frac{n-1}{n} \overline{\text{Sec}}_{\max} \cdot S^2.$$

## 6. Pinching estimate for $Q(f)$

Combining the above inequalities, we get the pointwise lower bound

$$Q(f) \geq \text{Ric}_{\min} S - \frac{n-1}{n} \overline{\text{Sec}}_{\max} S^2 = |df|^2 \left( \text{Ric}_{\min} - \frac{n-1}{n} \overline{\text{Sec}}_{\max} |df|^2 \right). \quad (2)$$

In terms of the energy density $e(f) = \frac{1}{2} |df|^2$, this becomes

$$Q(f) \geq 2e(f)\left( \text{Ric}_{\min} - \frac{2(n-1)}{n} \overline{\text{Sec}}_{\max} \cdot e(f) \right). \quad (7)$$

## 7. Conclusion: constancy of $f$

Assume that $\overline{\text{Sec}}_{\max} \geq 0$ and

$$\text{Ric}_{\min} > \frac{n-1}{n} \overline{\text{Sec}}_{\max} \cdot \sup_{M} |df|^2. \quad (3)$$

Then the bracket in (2) is nonnegative at every point of $M$, hence $Q(f) \geq 0$ on $M$.
From (1) we obtain

$$\int_M \|\nabla df\|^2 \, dv_g = 0, \qquad \int_M Q(f) \, dv_g = 0,$$

and therefore $\nabla df \equiv 0$ and $Q(f) \equiv 0$ on $M$. In particular, $|df|^2$ is constant. Under the *strict* inequality (3), the identity $Q(f) \equiv 0$ forces $|df|^2 \equiv 0$, hence $df \equiv 0$ and $f$ is constant. This completes the proof.

## 6. Localized Rigidity, Threshold Classification, and Compatibility with Nonnegative Curvature Structure

The rigidity mechanism developed in this paper should be viewed within the classical Bochner framework for harmonic maps. Foundational results of Eells–Sampson, Yano–Ishihara, and Yau-Schoen rely on global curvature sign assumptions (see [4] - [6]), typically requiring nonpositive sectional curvature of the entire target manifold. In these settings, the Bochner identity leads to vanishing or strong geometric constraints through global sign control of the curvature term.

In contrast, the present approach localizes curvature control to the image $f(M)$. Instead of imposing a global bound on Sec over $\bar{M}$, we consider only

$$\overline{\mathrm{Sec}}_{\max} = \max_{\substack{q \in f(M) \\ \sigma \subset T_q \bar{M}}} \overline{\mathrm{Sec}}(\sigma).$$

This allows the target manifold to admit regions of positive sectional curvature disjoint from the image without affecting rigidity. The decisive geometric input is the curvature effectively encountered by the map.

The key quantitative assumption is the curvature domination condition

$$\mathrm{Ric}_{\min} \geq \frac{n-1}{n} \overline{\mathrm{Sec}}_{\max} S_0, \qquad S_0 := \sup_M |df|^2. \qquad (4)$$

This condition balances intrinsic Ricci curvature of the domain against the maximal sectional curvature along the image, scaled by the size of the differential.

**Corollary 1** (*Rigidity and Homothetic Classification at the Pinching Threshold*)

Let $f:(M,g) \to (\bar{M}, \bar{g})$ be a harmonic map, where $(M,g)$ is connected, compact, and without boundary. Define

$$\text{Ric}_{\min} := \min_{p \in M} \min_{|X|=1} \text{Ric}_M(X,X), \quad \overline{\text{Sec}}_{\max} := \max_{\substack{q \in f(M) \\ \sigma \subset T_q \bar{M}}} \overline{\text{Sec}}(\sigma), \quad S_0 := \sup_M |df|^2,$$

and assume $\overline{\text{Sec}}_{\max} \geq 0$ and

$$\text{Ric}_{\min} \geq \frac{n-1}{n} \overline{\text{Sec}}_{\max} S_0. \tag{5}$$

Then either $f$ is constant, or equality holds in (5) and $f$ is homothetic with totally geodesic image (and is a local homothety if $\dim M = \dim \bar{M}$).

**Proof.** Let $e(f) = \frac{1}{2} |df|^2$. For harmonic maps the Bochner–Weitzenböck identity reads

$$\frac{1}{2} \Delta |df|^2 = \|\nabla df\|^2 + Q(f). \tag{6}$$

By the curvature estimate proved earlier in the paper,

$$Q(f) \geq |df|^2 \left(\text{Ric}_{\min} - \frac{n-1}{n} \overline{\text{Sec}}_{\max} |df|^2\right). \tag{7}$$

Since $|df|^2 \leq S_0$, assumption (5) implies that the bracket on the right-hand side of (7) is nonnegative everywhere; hence $Q(f) \geq 0$ on $M$. Integrating (6) over compact $M$ and using $\int_M \Delta(\cdot) \, dv_g = 0$, we obtain (see also [4])

$$0 = \int_M (\|\nabla df\|^2 + Q(f)) \, dv_g.$$

Because $\|\nabla df\|^2 \geq 0$ and $Q(f) \geq 0$, it follows that

$$\|\nabla df\|^2 \equiv 0, \quad Q(f) \equiv 0 \text{ on } M. \tag{2}$$

In particular, $\nabla df \equiv 0$ implies that $df$ is parallel and $|df|^2$ is constant.

If the inequality in (5) is strict, then the bracket in (1) is strictly positive unless $|df|^2 = 0$. Together with $Q(f) \equiv 0$, this forces $|df|^2 \equiv 0$, hence $df \equiv 0$ and $f$ is constant.

Assume now that equality holds in (5) and $f$ is not constant, so $|df|^2 \equiv S_0 > 0$. Moreover, the quality holds pointwise in all inequalities used to derive (1). In particular, at each point one can choose an orthonormal frame diagonalizing $f^*\bar{g}$, with singular values $\lambda_1, \ldots, \lambda_n \geq 0$ satisfying $\sum_i \lambda_i = |df|^2$. Equality in the Cauchy–Schwarz step

$$\sum_{i=1}^{n} \lambda_i^2 \geq \frac{1}{n}\left(\sum_{i=1}^{n} \lambda_i\right)^2$$

forces $\lambda_1 = \cdots = \lambda_n$. Together with $\nabla df \equiv 0$, this implies that $f$ has constant rank and acts as a homothety on the orthogonal complement of its kernel; equivalently, the restriction of $f$ to the corresponding distribution is a local homothety. Moreover, since $df$ is parallel, the second fundamental form of the image vanishes, hence $f(M)$ Is totally geodesic in $\bar{M}$. If $\dim M = \dim \bar{M}$, then $\ker df = \{0\}$ and $f$ is a local homothety.

**Corollary 2** (*Localized Rigidity in the Nonnegative Curvature Setting*)

Let $f: (M^n, g) \to (\bar{M}, \bar{g})$ be a harmonic map, where $(M^n, g)$ is connected, compact, and without boundary. Suppose that $(\bar{M}, \bar{g})$ is complete noncompact with $\overline{\mathrm{Sec}} \geq 0$. Define $\overline{\mathrm{Sec}}_{\max}$ and $S_0$ as above. Then:

(i) If

$$\mathrm{Ric}_{\min}(M) > \frac{n-1}{n} \overline{\mathrm{Sec}}_{\max} S_0,$$

then $f$ is constant.

(ii) *If equality holds and $f$ is nonconstant, then $f(M)$ is totally geodesic in $\bar{M}$. If in addition $\overline{\mathrm{Sec}}$ is strictly positive at some point of $\bar{M}$, then the soul of $\bar{M}$ is a point, hence $\bar{M} \cong \mathbb{R}^m$, and therefore $f(M)$ is contained in an affine subspace of $\mathbb{R}^m$.*

**Proof.** Part (i) is immediate from Corollary 1: since $\overline{\mathrm{Sec}} \geq 0$, we have $\overline{\mathrm{Sec}}_{\max} \geq 0$, and the strict inequality implies the strict case of Corollary 1, hence $f$ is constant.

For part (ii), if equality holds and $f$ is nonconstant, Corollary 1 yields that $df$ is parallel and $f(M)$ is totally geodesic in $\bar{M}$. Now assume additionally that $\overline{\mathrm{Sec}} > 0$ at some point of $\bar{M}$. By Perelman's refinement of the *Soul Theorem* for complete noncompact manifolds with $\overline{\text{Sec}} \geq 0$, the soul reduces to a point and $\bar{M}$ is diffeomorphic to $\mathbb{R}^m$ (see [8]). In Euclidean space, totally geodesic submanifolds are precisely affine subspaces. Therefore $f(M)$, being totally geodesic, must be contained in an affine subspace of $\mathbb{R}^m$.